\documentclass[12pt]{article}

\usepackage{amsmath}
\usepackage{amsfonts}

\def\mnalg{{\cal A}}
\def\sbalg{{\cal B}}

\def\Bbb{\mathbb}

\def\Cmpx{{\Bbb{C}}}

\def\Intg{{\Bbb{Z}}}
\def\spanone{{\Bbb{I}}}
\def\span{{\hbox{span}}}
\def\ad{\hbox{\rm ad}}
\def\tr{\hbox{\rm tr}}

\def\Der{{\hbox{\rm Der}}}
\def\rank{{\cal R}}
\def\Lie{{\cal L}}

\def\cnj{\overline}

\def\qed{\hfill \vrule height 0.7 em width 0.7 em}
\def\tfrac#1#2{{\textstyle{\frac{#1}{#2}}}}

\newtheorem{theorem}{Theorem}
\newtheorem{lemma}[theorem]{Lemma}

\numberwithin{equation}{section}

{%
\refstepcounter{theorem}%
{\vskip 1em\noindent \bf Example\ \thetheorem :}}%
{}

\newenvironment{proof}[1]%
{\vskip 1em\noindent {\bf Proof #1\par}}%
{\qed\vskip 1em}

\title{Non Commutative Differential Geometry, and the Matrix
Representations of Generalised Algebras}

\author{{\bf J. Gratus}%
\\
Laboratoire de Gravitation et Cosmologie Relativistes%
\thanks{\it Laboratoire associ\'e au CNRS {\rm URA 769}}
\\
Tour 22/12 4eme etage, Boite Courrier142, 4pl Jussieu. F75252
Paris
\\
email: gratus@ccr.jussieu.fr
}

\date{14 March 1997}

\begin{document}
\maketitle


\begin{abstract}

The underlying algebra for a noncommutative geometry is taken to be a
matrix algebra, and the set of derivatives the adjoint of a subset of
traceless matrices.  This is sufficient to calculate the dual 1-forms,
and show that the space of 1-forms is a free module over the algebra
of matrices. The concept of a generalised algebra is defined and it is
shown that this is required in order for the space of 2-forms to
exist. The exterior derivative is generalised for higher order forms
and these are also shown to be free modules over the matrix algebra.
Examples of mappings that preserve the differential structure are
given.  Also given are four examples of matrix generalised algebras,
and the corresponding noncommutative geometries, including the cases
where the generalised algebra corresponds to a representation of a Lie
algebra or a $q$-deformed algebra.

\end{abstract}

\begin{flushleft}
Key Words: \\
\qquad (1) Noncommutative Geometry \\
\qquad (2) Generalised Algebra \\[1em]
MCS 1991: \\
Primary: 
17A30 {General nonassociative rings: Algebras
satisfying 
other identities} \\
Secondary:  
81S05 {General quantum mechanics and problems of
quantization: 
Commutation relations} \\
\qquad\qquad
58A10 {General theory of differentiable manifolds:
Differential forms} \\[1em]
PACS: \\
Primary:  03.65.Fd Algebraic Methods (Quantum Mechanics) \\
Secondary:  04.60.Ne Lattice and discreet methods (Quantum Gravity) \\
\qquad\qquad 02.10.Sp Linear and multilinear algebra: matrix theory
(finite and infinite) \\
\end{flushleft}

\newpage

%

\section{Introduction}

To define a noncommutaive geometry or differential calculus, it is
first necessary to introduce an algebra $\mnalg$ that will replace the
algebra of functions. There is a unique universal differential
calculus for which all calculi are quotients.  There are several
methods of defining the quotient map necessary for the space of
1-froms. The method we use follows \cite{Madore_Dim1} by
constructing it with respect to a subspace $\sbalg$ of $\mnalg$.

In the early days \cite{Connes1,Connes2} $\sbalg$ was taken to be
$\mnalg$ itself.  Later, \cite[chapter 3]{Madore_book} examples where
$\sbalg$ formed a Lie algebra, or some other algebraic
relationship such as $[p,x]=1$ as in quantum mechanics, or
$xy=qyx$ as in $q$-deformed algebras, were studied.

For each subspace $\sbalg$ one could construct a co-frame. This
co-frame is loosely analogous to the orthonormal co-frame used in
normal differential geometry.  By quotienting the universal calculus
one could then construct the set of 2-forms and higher order forms.

It was discovered that if $\sbalg=\mnalg$, or $\sbalg$ formed a Lie
algebra or a quantum algebra then one could consistently impose the
condition that the co-frame basis elements of the exterior algebra
anticommute. Whilst for $q$-deformed algebras the basis elements of
the exterior algebra $q$-anticommute.

It is only recently \cite{Madore_Dim1,Madore_Dim2} that people have
looked at a general $\sbalg$. They showed that in order for forms of
order 2 and above to exist puts constraints on the elements of
$\sbalg$.  However these constraints have not been pursued.

In this paper we impose the condition $\mnalg=M_m(\Cmpx)$ and that all
the elements in $\sbalg$ are traceless. In section \ref{ch_1F} we show
that this is a sufficient condition for the co-frame to exist.
However this condition is not a necessary condition for the co-frame
to exist. To show this we give some examples where $\mnalg\ne
M_n(\Cmpx)$, some of which have a co-frame and others which do not.
Since $M_m(\Cmpx)$ is a finite approximation to the infinite
dimensional space of functions it is hoped that this procedure can be
used as an alternative to the theory of renormalisation or lattice
QFT.

In section \ref{ch_GA} we introduce the concept of a ``generalised
algebra''.  This is an algebraic structure that includes commutative
algebras, anti-commutative algebras, Lie algebras, Clifford algebras
and $q$-deformed algebras as examples. Each generalised algebra has a
specific rank and the space of 2-forms is a free module over $\mnalg$
of rank equal to the rank of the generalised algebra.  In section
\ref{ch_2F} we show that for 2-forms and higher forms to exist
$\sbalg$ must form a generalised algebra.  In section \ref{ch_HO}, we
then give the structure of the higher order forms, all of which are
also free modules over $\mnalg$, and an explicit expression for the
exterior derivative. In section \ref{ch_map} we give a couple of
simple examples of maps between generalised algebras which are
$d$-homomorphism, i.e. they preserve the differentiable stucture.

To elucidate the relationship between the generalised algebra of
$\sbalg$ and the space of 2-forms we give, in section \ref{ch_eg},
four examples: Much emphasis has been place on the case that $\sbalg$
form a Lie algebra.  Especially since $su(2)$ corresponds to the fuzzy
or non commutative version of the sphere \cite[chapter
7.2]{Madore_book} and $su(4)$ is an analogue of the Euclidianised
compactified Minkowski space \cite{Grosse1}.  Another example is that
of the $q$-deformed algebra, this has a finite dimensional
representation only if there exists an $m\in\Intg$ such that $q^m=1$.
Finally a $\sbalg$ is given of dimension 3, and rank 1 which may be
thought of as the fuzzy ellipse.

For further references, and history of this subject the reader is
asked to read the book \cite{Madore_book}.


\subsection{Note on Notation}

Unless otherwise stated $\mnalg=M_m(\Cmpx)$. $\sbalg\subset\mnalg$ is
a subspace of dimension $n$ of traceless matrices and $\lambda_a$ is
a basis for $\sbalg$.  Early Roman letters used as indices
$a,b,\ldots$ run over $1\ldots n$, and we use the Einstein summation
convention so that the is implicit summation if one index is high and
the other low.  The indices $r,s=1\ldots \rank$, while Greek
indices $\mu,\nu=1\ldots m^2$ and also follow the summation
convention.


\section{Generalised Algebras}

\label{ch_GA}

Given an Algebra $\mnalg$ with a unit, a subspace of that algebra
$\sbalg\subset\mnalg$ of finite dimension $n$ is said to be a {\bf
Generalised Algebra} of rank $\rank$ if for any basis
$\{\lambda_a\}_{a=1\ldots n}$ of $\sbalg$, there exist a
$n^2\times\rank$ matrix of rank $\rank$ given by $(\alpha^{ab}_r)$
such that
\begin{eqnarray}
\alpha^{ab}_r \lambda_a\lambda_b \in\sbalg\oplus\spanone 
\label{GA_alpha_sb}
\end{eqnarray}
where $\spanone=\span \{1\}$. Here $a,b$ are summed over $1\ldots n$,
$r=1\ldots\rank$, and $\rank\le n^2$.  As stated in the introduction
throughout this article we shall assume that $\mnalg=M_m(\Cmpx)$.  We
can think of (\ref{GA_alpha_sb}) simply as a set of relationships on
the independent matrices $\{\lambda_a\}$.  Alternatively we can think
of a generalised algebra as an abstract vector space, with the only
products defined being those defined by (\ref{GA_alpha_sb}). The
mapping that takes the elements of the generalised algebra into
matrices can be thought of as a (matrix) representation of the
underlying generalised algebra. In the same way as we think of Lie
algebras and Clifford algebras as being the fundamental object, and
the $\gamma$ matrices as merely a representation.

The matrix algebra $M_m(\Cmpx)$ comes equipped with an inner product
\begin{eqnarray}
\langle f,g \rangle &=& \tr(f^\dagger g)
\label{GA_inner_product}
\end{eqnarray}
Since trace is defined we shall assume that all elements in $\sbalg$
are traceless matrices.
Since the inner product is positive definite its restriction onto
$\sbalg$ is also positive definite and so the matrix
\begin{eqnarray}
g_{ab}=\langle\lambda_a,\lambda_b\rangle
\end{eqnarray}
is positive definite and Hermitian, $g_{ab}=\cnj{g_{ba}}$.  We label
its inverse by $g^{ab}$, and define the elements
$\{\lambda^a\in\sbalg\}$ dual to $\{\lambda_b\}$ by
\begin{eqnarray}
\lambda^a = g^{ba} \lambda_b
\qquad\hbox{so}\qquad
\langle\lambda^a,\lambda_b\rangle=\delta^a_b
\end{eqnarray}

It is also useful to define the orthogonal projections onto, and
perpendicular to $\sbalg$
\begin{eqnarray}
\eta &:& \mnalg\mapsto\sbalg\subset\mnalg\cr
\eta(f) &=& \langle \lambda^a , f \rangle \lambda_a \\
\eta^\perp &:&  \mnalg\mapsto\mnalg\cr
\eta^\perp(f) &=& f - \eta(f) - \tfrac1m\tr(f)
\end{eqnarray}
By taking the trace of (\ref{GA_alpha_sb}) and its orthogonal
projection onto $\sbalg$, we get
\begin{eqnarray}
&&\alpha^{ab}_r \left(\lambda_a\lambda_b - 
\langle\lambda^c,\lambda_a\lambda_b \rangle\lambda_c
- \tfrac1m \tr(\lambda_a\lambda_b) \right) = 0
\label{GA_alpha}
\\
&&\alpha^{ab}_r \eta^\perp(\lambda_a\lambda_b) = 0
\label{GA_eta}
\end{eqnarray}

We shall see in section \ref{ch_2F} that it is useful to construct
the $n^2\times n^2$ projection matrix $P^{ab}{}_{cd}$ of rank $\rank$
so that we can write (\ref{GA_alpha_sb}) as
\begin{eqnarray}
P^{ab}{}_{cd}\left(\lambda_a\lambda_b - 
\langle\lambda^e,\lambda_a\lambda_b \rangle\lambda_e
- \tfrac1m \tr(\lambda_a\lambda_b) \right) = 0
\label{GA_P}
\end{eqnarray}
For this we simply require a $n^2\times\rank$ matrix $\beta^r_{ab}$
such that
\begin{eqnarray}
\beta^r_{ab} \alpha^{ab}_s &=& \delta^r_s \qquad r,s=1\ldots \rank
\label{GA_def_beta}
\\
P^{ab}{}_{cd} &=& \sum_{r=1}^{\rank} \alpha^{ab}_r \beta^r_{cd} 
\label{GA_def_P_albe}
\end{eqnarray}
The choice of $\sbalg$ constrains, but does not completely determine
$\beta^r_{cd}$ and thus $P^{ab}{}_{cd}$.  From (\ref{GA_def_beta}) we
see that there are $n^2(n^2-\rank)+\rank^2$ linear constraints on the
$n^4$ elements of $P^{ab}{}_{cd}$.

We note that for the given $\sbalg$, $\rank$ might not be maximal,
i.e. there may exist other independent equations of the form
(\ref{GA_alpha}) which we have chosen to ignore. Therefore we
have the inequality
\begin{eqnarray}
\dim\left(\span \Big\{
\{\lambda_a\lambda_b\}_{a,b=1\ldots n},
\{\lambda_a\}_{a=1\ldots n},
1\Big\}\right)
\le
n^2 + n + 1 - \rank
\end{eqnarray}


\section{The Differential Calculi: 1 forms}
\label{ch_1F}

Let $\mnalg$ be any unital associative $\star$-algebra.  Of the many
differential calculi which can be constructed over $\mnalg$ the
largest is the differential envelope or universal differential
calculus $(\Omega^\star_u(\mnalg), d_u)$. Every other differential
calculus can be considered as a quotient of it.  For the definitions
refer to, for example, \cite{Connes1,Connes2}\cite[chapter
6.1]{Madore_book}. Let $(\Omega^\star(\mnalg), d)$ be another
differential calculus over $\mnalg$. Then there exists a unique
surjective $d_u$-homomorphism $\phi$
\begin{eqnarray}
\def\normalbaselines{\baselineskip=18pt}
\begin{matrix}
\mnalg &\buildrel d_u \over \longrightarrow &\Omega_u^1(\mnalg)
         &\buildrel d_u \over \longrightarrow &\Omega_u^2(\mnalg)
         &\buildrel d_u \over \longrightarrow &\cdots            
\cr
\parallel&&\phi_1 \downarrow \phantom{\phi_1}
         &&\phi_2 \downarrow \phantom{\phi_2}                      
\cr
\mnalg &\buildrel d \over \longrightarrow &\Omega^1(\mnalg)
         &\buildrel d \over \longrightarrow &\Omega^2(\mnalg)
         &\buildrel d \over \longrightarrow &\cdots
\end{matrix}
\def\normalbaselines{\baselineskip=12pt}                          
\label{1F_comut_diag} 
\end{eqnarray}
of $\Omega^\star_u(\mnalg)$ onto $\Omega^\star(\mnalg)$. It is given by
\begin{eqnarray}
\phi (d_u \xi) =  d \xi.
\end{eqnarray}
The restriction $\phi_p$ of $\phi$ to each $\Omega^p_u$ is defined by
\begin{eqnarray}
\phi_p(f_0 d_u f_1 \cdots d_u f_p) = f_0 df_1 \cdots df_p.
\end{eqnarray}
making $\Omega^\star(\mnalg)$ a bimodule over $\mnalg$.

Let us define $\Omega^1_\sbalg=\Omega^1_\sbalg(\mnalg)$, with respect
to $\sbalg\subset\mnalg$ by requiring
\begin{eqnarray}
\ker(\phi_1) &=& 
\left\{\sum_i f_i\;dg_i,\hbox{ with $f_i,g_i\in\mnalg$}
\ \bigg|\ \sum_i f_i[h,g_i] = 0\ \forall h\in\sbalg
\right\}
\label{1F_def_kerphi}
\end{eqnarray}
This is sufficient to define $\Omega^1_\sbalg$. 
We define the set of derivations 
\begin{eqnarray}
\Der_\sbalg &=& \{\ad(h)\ |\ h\in\sbalg\}   
\end{eqnarray}
this is a
complex vector space of dimension $n$. 
We now have the contraction given by
\begin{eqnarray}
\cdot:\Omega^1_\sbalg\times\Der_\sbalg &\mapsto& \mnalg
\cr
\sum_i f_i\;dg_i \cdot \ad(h) &=& \sum_i f_i[h,g_i] 
\end{eqnarray}
Which satisfies
\begin{eqnarray}
(f \xi g) \cdot X = f (\xi\cdot X) g
\qquad\forall f,g\in\mnalg,\ \xi\in\Omega^1_\sbalg,\ X\in\Der_\sbalg
\label{1F_def_dot}
\end{eqnarray}


From (\ref{1F_def_kerphi}) we see that for
$\xi\in\Omega^1_\sbalg$ then 
$\xi\cdot X=0$ for all $X\in\Der_\sbalg$ implies $\xi=0$.
Thus there is an injective linear map from $\Omega^1_\sbalg$ into the
dual over $\mnalg$ of $\Der_\sbalg$:
\begin{eqnarray}
\Omega^1_\sbalg
\hookrightarrow
\Der^*_\sbalg 
\;{\buildrel \hbox{\scriptsize def}\over =}\;
\{ \xi:\Der_\sbalg\mapsto\mnalg\ |\ \xi\hbox{ is linear } \}
\end{eqnarray}
We say that $\Omega^1_\sbalg$ {\bf has a co-frame} if
$\Omega^1_\sbalg=\Der^*_\sbalg$. 

Given the basis $\{\lambda_a\}_{a=1\ldots n}$ of $\sbalg$, we have the
basis $\{e_a=\ad(\lambda_a)\}_{a=1\ldots n}$ of $\Der_\sbalg$.
If $\Omega^1_\sbalg$ has a co-frame then we can define the
the {\bf co-frame} forms $\theta^b$ to be dual to $e_a$ by
\begin{eqnarray}
\theta^b\cdot e_a &=& \delta^a_b
\label{1F_def_theta_a}
\end{eqnarray}
From (\ref{1F_def_theta_a}) and (\ref{1F_def_dot}) we have
\begin{eqnarray}
\theta^a f &=& f \theta^a \qquad\forall f\in\mnalg
\end{eqnarray}
We define the form $\theta$ to be
\begin{eqnarray}
\theta\cdot\ad(h) &=& -h \qquad\forall h\in\sbalg
\end{eqnarray}
which has the following identities
\begin{eqnarray}
-[\theta,f] &=& d f \qquad\forall f\in\mnalg \\
\theta &=& -\lambda_a\theta^a
\end{eqnarray}
The relationship between these objects and those found in normal
differential geometry are vague.  The derivations $\{e_i\}$ are said
to be analogous to the orthonormal frame for normal differential
geometry whilst $\{\theta^a\}$ correspond to its dual co-frame.  There
is no analogy to the form $\theta$.


As already stated, in this article we shall take $\mnalg=M_m(\Cmpx)$
and $\sbalg\subset\mnalg$ as an $n$ dimensional subspace of traceless
matrices. This is because of:

\begin{theorem}
\label{1F_thm_frame}
Assuming $\mnalg=M_m(\Cmpx)$ and $\sbalg\subset\mnalg$ is
an $n$ dimensional subspace of traceless matrices then:

\begin{list}{$\bullet$}{}

\item
There are exact expressions for $\theta^a$ and $\theta$ given by
\begin{eqnarray}
\theta^a  &=&
\gamma_\nu \lambda^{a\dagger} d\gamma^{\nu\dagger} 
\label{1F_theta_a}
\\
\theta &=&
 \tfrac1{m}\gamma_\mu d\gamma^{\mu\dagger}
= -\tfrac1{m}d\gamma^{\mu\dagger} \gamma_\mu
\label{1F_theta_sum_dlamda}
\end{eqnarray}
where $\{{\gamma_\mu}\}_{\mu=1\ldots m^2}$ refer to any
basis of $\mnalg=M_m(\Cmpx)$, and we set
$\{\gamma^\nu\}_{\mu=1\ldots m^2}$ to be its dual, so that
$\langle \gamma^\nu , \gamma_\mu \rangle= \delta^\nu_\mu$.

\item
$\Omega^1_\sbalg$ has a co-frame. 

\item
$\Omega^1_\sbalg$ is a free module of rank $n$ over $\mnalg$. viz
\begin{eqnarray}
\Omega^1_\sbalg&= \otimes^n\mnalg
\label{1F_free_mod}
\end{eqnarray}
\end{list}

\end{theorem}


Before proving these we observe the following lemma

\begin{lemma}
\begin{eqnarray}
&&\gamma_\mu f \gamma^{\mu\dagger} = \tr(f) 
\qquad\forall f\in\mnalg
\label{1F_lm_tr_eq1}
\\
&&[ f , \gamma_\mu g \otimes \gamma^{\mu\dagger}] = 0
\qquad\forall f,g\in\mnalg
\label{1F_lm_eq2}
\end{eqnarray}
where $\{\gamma_\mu,\gamma^\nu\}_{\mu=1\ldots m^2}$ 
as in theorem \ref{1F_thm_frame}.
\label{1F_lm_tr}
\end{lemma}

\begin{proof}{of lemma \ref{1F_lm_tr}}
First note that these are independent of the choice of
basis $\gamma_\mu$. 

Now choose the basis $\{E_{ij}\}_{i,j=1\ldots m}$ to be the matrix
with a 1 in the $i$th row and the $j$th column, the natural basis for
the $m\times m$ matrices, so that $E_{ij}E_{kl}=E_{il}\delta_{jk}$.
These elements are orthonormal with respect to the trace inner product
so $E_{ij}$ is dual to itself.  (During this proof indices $i,j,k,l$
are summed from $1\ldots m$.)  Now let $f=f_{kl}E_{kl}$ with
$f_{kl}\in\Cmpx$ so
\begin{eqnarray}
E_{ij}f E_{ji} &=& 
E_{ij}f_{kl}E_{kl} E_{ji} 
=
f_{kl} \delta_{jk}\delta_{jl} E_{ii} 
=
f_{jj} E_{ii} \cr
&=&
\tr(f)
\end{eqnarray}
Whilst
\begin{eqnarray}
f E_{ij} g \otimes E_{ji} - E_{ij} g \otimes E_{ji} f
&=&
f_{kl}E_{kl} E_{ij} g \otimes E_{ji} - E_{ij} g \otimes E_{ji}
f_{kl}E_{kl} \cr
&=&
f_{kl}(E_{kj} \delta_{li} g \otimes E_{ji} -
E_{ij} g \otimes E_{jl}\delta_{ik} ) \cr
&=&
f_{kl}(E_{kj} g \otimes E_{jl} -
E_{kj} g \otimes E_{jl} ) \cr
&=&
0
\end{eqnarray}
\end{proof}


\begin{proof}{of theorem \ref{1F_thm_frame}}

From (\ref{1F_lm_tr_eq1}) we have
\begin{eqnarray}
\gamma_\nu \lambda^{a\dagger} {d}\gamma^{\nu\dagger} 
&=&
\gamma_\nu \lambda^{a\dagger} 
(\gamma^{\nu\dagger}\theta - \theta\gamma^{\nu\dagger})
\cr
&=&
(-\gamma_\nu\lambda^{a\dagger}\gamma^{\nu\dagger}
\lambda_b
+
\gamma_\nu\lambda^{a\dagger}
\lambda_b\gamma^{\nu\dagger})\theta^b
\cr
&=&
(-\tr(\lambda^{a\dagger})\lambda_b + \tr(\lambda^{a\dagger}\lambda_b)
\theta^b
\cr
&=&
\langle\lambda^a,\lambda_b\rangle
\theta^b
=
\delta^a_b\theta^b
=
\theta^a
\end{eqnarray}
hence (\ref{1F_theta_a}). Also
\begin{eqnarray}
\gamma_\mu {d}\gamma^{\mu\dagger}
&=&
\gamma_\mu [\lambda_a,\gamma^{\mu\dagger}] \theta^a 
\cr
&=&
-m \lambda_a \theta^a  = m \theta
\end{eqnarray}
whilst
\begin{eqnarray}
0 &=& {d}(\gamma_\mu \gamma^{\mu\dagger})
= \gamma_\mu {d}\gamma^{\mu\dagger} +
{d}(\gamma_\mu)\gamma^{\mu\dagger}
\end{eqnarray}
Thus (\ref{1F_theta_sum_dlamda}). Given any linear map
$\xi:\Omega^1_\sbalg\mapsto\mnalg$ then let $\xi_a=\xi\cdot
e_a\in\mnalg$ then $\xi=\xi_a\theta^a$ so $\Omega^1_\sbalg$ has a
co-frame, and is also a free module over $\mnalg$ with rank $n$ and
basis $\{\theta^a\}$.

\end{proof}

The elements of $\Omega_u^1(\mnalg)$ which map onto
$\theta^a$ and $\theta$ by the projection $\phi^1$ are given by:
\begin{eqnarray}
\theta^a_u &=& \gamma_\mu\lambda^{a\dagger}\otimes \gamma^{\mu\dagger}
=
\gamma_\mu\lambda^{a\dagger} d_u\gamma^{\mu\dagger}
\label{1F_theta_u_a}
\\
\theta_u &=&  \tfrac1m \gamma_\mu \otimes \gamma^{\mu\dagger} 
- 1\otimes 1 
=
\tfrac1m \gamma_\mu d_u\gamma^{\mu\dagger}
\label{1F_theta_u}
\\
\noalign{\hbox{so}}
\phi^1(\theta^a_u) &=& \theta^a
\qquad\hbox{and}\qquad
\phi^1(\theta_u) = \theta
\end{eqnarray}
As $\phi$ is not injective, $\theta^a_u$ and $\theta_u$ are not
unique. However $\theta^a_u=\eta^\star(\theta^a)$ (see example in
section \ref{ch_eg_A0}), and
$\theta_u$ does satisfy
\begin{eqnarray}
-[\theta_u,f] &=& d_u f
\label{1F_comm_theta_u}
\end{eqnarray}
which is shown by using (\ref{1F_lm_eq2}).


\vskip 1em
{\large\bf Counter Examples}
\label{1F_counter_eg}

If one does not require both that $\mnalg=M_m(\Cmpx)$ and that all the
elements in $\sbalg$ are traceless then the question of whether
$\Omega^1_\sbalg$ has a co-frame is non trivial. 
Here are some examples where 
$\Omega^1_\sbalg$ does not have a co-frame:
\begin{list}{$\bullet$}{}

\item
$\mnalg$ is Abelien.

\item
$\mnalg=M_m(\Cmpx)$ but $1\in\sbalg$. This is because $\ad(1)=0$.

\item
$\mnalg=\{\hbox{space of operators generated by $x$ and $p$ where
$[p,x]=1$}\}$ and $\sbalg=\span\{p,p^2,x\}$.

\end{list}

Whilst on the contrary $\Omega^1_\sbalg$ does have a co-frame
\begin{list}{$\bullet$}{}

\item 
$\mnalg=M_m(\Cmpx)$ but $\sbalg=\span\{1+x,y,z\}$, where $\{x,y,z\}$
is a representation of $su(2)$.

\item 
$\mnalg=\{\hbox{space of operators generated by $x$ and $p$ where
$[p,x]=1$}\}$ and $\sbalg=\span\{p,x\}$. This is the Heisenberg quantum
algebra.

\end{list}


\section{$\Omega^2_\sbalg$ and Generalised Algebras}
\label{ch_2F}

Having constructed the set of 1-forms $\Omega^1_\sbalg$ we turn our
attention to $\Omega^2_\sbalg$, the structure of which is
given by the following theorem:

\begin{theorem}
\label{2F_thm_dtheta}
Given $d:\Omega^1_\sbalg\mapsto\Omega^2_\sbalg$ obeys
(\ref{1F_comut_diag}) we have the following
\begin{eqnarray}
&&{d}\theta + \theta{}^2 = - \tfrac1m \tr(\lambda_a\lambda_b)
\theta^a \theta^b
\label{2F_dtheta_theta_sqr}
\\
&&{d}\theta^a  =
-[\theta,\theta^a ]
-\langle\lambda^a ,\lambda_b \lambda_c\rangle
\theta^b \theta^{c}
\label{2F_dtheta_a}
\\
&&\eta^\perp(\lambda_a\lambda_b)\theta^a \theta^b 
= 0
\label{2F_eta_theta2}
\end{eqnarray}
where $[\bullet,\bullet]$ is the graded commutator.  If the
contraction of two-forms on pairs of vectors obeys the 2-form version
of (\ref{1F_def_dot})
\begin{eqnarray}
&&(f \xi g)\cdot (X,Y) = f (\xi\cdot(X,Y)) g
\qquad\forall f,g\in\mnalg,\ \xi\in\Omega^2_\sbalg,\ X,Y\in\Der_\sbalg
\label{2F_def_dot}
\end{eqnarray}
then either $\dim(\Omega^2_\sbalg)=0$ or $\sbalg$ is a generalised
algebra. Let
\begin{eqnarray}
\theta^a \theta^b  \cdot 
({e}_c,{e}_d) &=&
P^{ab}{}_{cd}
\label{2F_def_Pabcd}
\end{eqnarray}
Viewing $P$ as an $n^2\times n^2$ matrix, if $P$ has rank $\rank$ then
$\Omega^2_\sbalg$ is a free module over $\mnalg$ of rank $\rank$. i.e.
\begin{eqnarray}
\Omega^2_\sbalg&=& \otimes^\rank \mnalg
\end{eqnarray}
\end{theorem}


\begin{proof}{} 
From (\ref{1F_comut_diag}) we have the standard relations on $d$
given by
\begin{eqnarray}
{d}({d}f) = 0 
\qquad\hbox{and}\qquad
{d}(f\;{d}h) = {d}f\;{d}h \
\qquad \forall f,h\in\mnalg
\label{2F_ddf_dfdg}
\end{eqnarray}
Using (\ref{1F_theta_sum_dlamda}) we have
\begin{eqnarray}
{d}\theta &=&
\tfrac1{m} {d}\gamma_\mu {d}\gamma^{\mu\dagger}
=
\tfrac1{m}
[\lambda_a,\gamma_\mu]
[\lambda_b,\gamma^{\mu\dagger}]
\theta^a \theta^b 
\cr
&=&
\tfrac1m( 
\gamma_\mu\lambda_a\gamma^{\mu\dagger}\lambda_b
-
\lambda_a\gamma_\mu\gamma^{\mu\dagger}\lambda_b
-
\gamma_\mu\lambda_a\lambda_b\gamma^{\mu\dagger}
+
\lambda_a\gamma_\mu\lambda_b\gamma^{\mu\dagger}
)
\cr
&=&
(-\lambda_a\lambda_b - \tfrac1m\tr(\lambda_a\lambda_b))
\theta^a \theta^b 
\end{eqnarray}
Hence (\ref{2F_dtheta_theta_sqr}).
From (\ref{1F_theta_a}) we have
\begin{eqnarray}
{d}\theta^a  &=&
{d}(\gamma_\nu \lambda^{a\dagger}  {d}\gamma^{\nu\dagger})
=
{d}(\gamma_\nu \lambda^{a\dagger} ) {d}\gamma^{\nu\dagger}
\cr
&=&
[\theta,\gamma_\nu \lambda^{a\dagger} ][\theta,\gamma^{\nu\dagger}]
=
[\lambda_b,\gamma_\nu \lambda^{a\dagger} ]
[\lambda_c,\gamma^{\nu\dagger}]\theta^b \theta^{c}
\cr
&=&
(
\lambda_b\gamma_\nu \lambda^{a\dagger}  \lambda_c\gamma^{\nu\dagger} 
-
\gamma_\nu\lambda^{a\dagger}  \lambda_b \lambda_c\gamma^{\nu\dagger}
-
\lambda_b\gamma_\nu \lambda^{a\dagger}  \gamma^{\nu\dagger}\lambda_c
+
\gamma_\nu\lambda^{a\dagger}  \lambda_b \gamma^{\nu\dagger}\lambda_c
)
\theta^b \theta^{c}
\cr
&=&
\Big(\lambda_b\tr(\lambda^{a\dagger} \lambda_c)
-\tr(\lambda^{a\dagger}  \lambda_b \lambda_c)
+\tr(\lambda^{a\dagger}  \lambda_b)\lambda_c
\Big)\theta^b \theta^{c}
\cr
&=&
\Big(\lambda_b\delta^a_c + \lambda_c \delta^a_b
-\langle\lambda^a ,\lambda_b \lambda_c\rangle
\Big)\theta^b \theta^{c}
\end{eqnarray}
Hence (\ref{2F_dtheta_a}).
From (\ref{2F_dtheta_theta_sqr}) we have
\begin{eqnarray}
{d}\theta &=&
-{d}(\lambda_a\theta^a )
=
-{d}\lambda_a\theta^a  - \lambda_a {d}\theta^a 
\cr
&=&
[\theta,\lambda_a]\theta^a 
+ \lambda_a 
\left(\theta\theta^a  
+\theta^a \theta 
+\langle\lambda^a,\lambda_b \lambda_c\rangle
\theta^b \theta^{c}
\right)
\cr
&=&
- 2 \theta^{2} + \eta(\lambda_b \lambda_c)\theta^b \theta^{c}
\end{eqnarray}
Comparing this with (\ref{2F_dtheta_theta_sqr}) gives
\begin{eqnarray}
- \theta^2 - \tfrac1m \tr(\lambda_a\lambda_b)
&=&
- 2 \theta^{2} + \eta(\lambda_b \lambda_c)\theta^b \theta^{c}
\end{eqnarray}
Hence (\ref{2F_eta_theta2})
From (\ref{2F_def_dot}) we have
\begin{eqnarray}
P^{ab}{}_{cd} f 
=
\theta^a \theta^b \cdot 
({e}_c,{e}_d) f  
= 
\theta^a \theta^b f \cdot 
({e}_c,{e}_d) 
&=& 
f\theta^a \theta^b  \cdot 
({e}_c,{e}_d) 
=
f P^{ab}{}_{cd}
\qquad\forall f\in\mnalg
\end{eqnarray}
hence $P^{ab}{}_{cd}$ is in the center of $\mnalg$ so is a multiple of
the unit element. 
Contracting this with (\ref{2F_eta_theta2}) gives
\begin{eqnarray}
\eta^\perp(P^{ab}{}_{cd}\lambda_a\lambda_b) &=& 0
\qquad\hbox{with }P^{ab}{}_{cd}\in\Cmpx 
\end{eqnarray}
Thus either $\dim(\Omega^2_\sbalg)=0$ or $\sbalg$ is a generalised
algebra. Any element $\xi\in\Omega^2_\sbalg$ can be written as
\begin{eqnarray}
\xi &=& \sum_\alpha f_0^{(\alpha)}\;df_1^{(\alpha)}\;df_2^{(\alpha)} 
\qquad\hbox{where}\
f_0^{(\alpha)},f_1^{(\alpha)},f_2^{(\alpha)}\in\mnalg
\cr
&=&
\sum_\alpha
f_0^{(\alpha)}[\lambda_a,f_1^{(\alpha)}][\lambda_b,f_2^{(\alpha)}]
\theta^a\theta^b
\end{eqnarray}
Thus $\Omega^2_\sbalg$ is a free module of $\mnalg$, with a basis
$\theta^a\theta^b$. From (\ref{2F_def_Pabcd}) we see the number of
independent sets of $\theta^a\theta^b$ is $\rank$.

\end{proof}

In order to be consistent with section \ref{ch_GA}
we shall assume that $P^2=P$ thus
\begin{eqnarray}
P^{ab}{}_{cd} \theta^c\theta^d =  \theta^a\theta^b
\end{eqnarray}
This is a special case of the results found in \cite{Madore_Dim1}
where $F^a{}_{bc} = \tfrac12 P^{de}{}_{bc}
\langle\lambda^a,\lambda_d\lambda_e\rangle$ and $K_{bc}=\tfrac1{2m}
P^{de}{}_{bc} \tr(\lambda_d\lambda_e)$


\section{Higher Order Forms}
\label{ch_HO}

The higher forms are still free modules over $\mnalg$ with the 
basis of $\Omega^p_\sbalg$ being a quotient of the set
\begin{eqnarray}
\{\theta^{a_1}\ldots\theta^{a_p}\}_{a_1\ldots a_p=1\ldots n}
\end{eqnarray}
The quotient being given by the extension of (\ref{2F_eta_theta2})
that adjacent pairwise contractions must vanish, viz
\begin{eqnarray}
\eta^\perp(\lambda_{a_q}\lambda_{a_{q+1}})
\theta^{a_1}\ldots\theta^{a_{q}}
\theta^{a_{q+1}}\ldots\theta^{a_r}
&=& 0
\qquad\forall q=1\ldots p-1
\label{HO_contract}
\end{eqnarray}
For this to give a non trivial free module requires that there exist
the $n^{p-2}\rank(p-1)$ complex numbers
\begin{eqnarray}
\Big\{
\varepsilon^{a_1\ldots a_{p-2}}_{r t}
\in\Cmpx\Big\}
\qquad\hbox{where}\ 
a_1\ldots a_{p-2}=1\ldots n,\ r=1\ldots\rank,\
t=1\ldots p-1
\end{eqnarray}
such that
\begin{eqnarray}
\sum_{r=1}^\rank
\alpha^{a_1 a_2}_r \varepsilon^{a_3\ldots a_p}_{r 1}
=
\sum_{r=1}^\rank
\alpha^{a_2 a_3}_r \varepsilon^{a_1 a_4\ldots a_p}_{r 2}
=
\ldots\ldots
=&&
\sum_{r=1}^\rank
\alpha^{a_{p-1} a_p}_r \varepsilon^{a_1\ldots a_{p-2}}_{r (p-1)}
\\
&&\forall
a_1\ldots a_p=1\ldots n,\ r=1\ldots\rank
\nonumber
\end{eqnarray}
The existence or otherwise of these $\varepsilon$'s and hence the rank
of $\Omega^p_\sbalg$ depends on the nature of $\alpha^{ab}_r$.
However if we do have a non trivial $\Omega^p_\sbalg$ then we have the
extension of $d$ given by the following theorem:

\begin{theorem}
\label{HO_thm}
Given that (\ref{HO_contract}) holds,
the extension to $d$ is given by
\begin{eqnarray}
{d} &:& \Omega^\star_\sbalg \mapsto \Omega^\star_\sbalg \cr
{d} &:& \Omega^p_\sbalg \mapsto \Omega^{p+1}_\sbalg \cr
{d}\xi &=& -[\theta,\xi] + \chi(\xi)
\label{HO_d}
\end{eqnarray}
where $[\bullet,\bullet]$ is the graded commutator and 
\begin{eqnarray}
\chi &:& \Omega^\star_\sbalg \mapsto \Omega^\star_\sbalg \cr
\chi &:& \Omega^p_\sbalg \mapsto \Omega^{p+1}_\sbalg \cr
\chi(f\theta^{a_1}\ldots\theta^{a_r})
&=&
 f\sum_{q=1}^p (-1)^{s+1}
\langle\lambda^{a_q},\lambda_b\lambda_c\rangle
\theta^{a_1}\ldots\theta^{a_{q-1}}
\theta^b\theta^c
\theta^{a_{q+1}}\ldots\theta^{a_p}
\end{eqnarray}
Both $d$ and $\chi$ are well defined and obey the graded Leibniz rule,
i.e:
\begin{eqnarray}
d(\zeta\xi)&=&
d(\zeta)\xi + (-1)^p\zeta d(\xi)
\qquad\forall \zeta\in\Omega^p_\sbalg,\, \xi\in\Omega^\star_\sbalg
\end{eqnarray}
$d$ obeys (\ref{1F_comut_diag}) and $\chi$ is left and right $\mnalg$
linear.
\end{theorem}


\begin{proof}{}
Now to show that they are well defined we note
\begin{eqnarray}
d(\eta^\perp(\lambda_a\lambda_b)\theta^a\theta^b) = 0
\end{eqnarray}
since
\begin{eqnarray}
\lefteqn{\chi(\eta^\perp(\lambda_a\lambda_b)\theta^a\theta^b)}\qquad
&&
\cr
&=&
\eta^\perp(\lambda_a\lambda_b)
\left( 
\langle \lambda^a,\lambda_c\lambda_d \rangle \theta^c\theta^d\theta^b
-
\langle \lambda^b,\lambda_c\lambda_d \rangle \theta^a\theta^c\theta^d
\right)
\cr
&=&
\eta^\perp\left(
\eta(\lambda_c\lambda_d)\lambda_b
-
\lambda_c\eta(\lambda_d\lambda_b)
\right)\theta^c\theta^d\theta^b
\cr
&=&
\eta^\perp\left(
\lambda_c\lambda_d\lambda_b -
\eta^\perp(\lambda_c\lambda_d)\lambda_b -
\tr(\lambda_c\lambda_d)\lambda_b
-
\lambda_c\lambda_d\lambda_b +
\lambda_c\eta^\perp(\lambda_d\lambda_b) +
\lambda_c\tr(\lambda_d\lambda_b)
\right)\theta^c\theta^d\theta^b
\cr
&=& 0
\end{eqnarray}
whilst
\begin{eqnarray}
[\theta,\eta^\perp(\lambda_a\lambda_b)\theta^a\theta^b]
=0
\end{eqnarray}
thus $d$ is well defined on 2-froms. Higher forms follow from
graded Leibniz. Also from graded Leibniz we have
\begin{eqnarray}
{d}(f_0{d}f_1\ldots {d}f_r)
&=&
{d}f_0{d}f_1\ldots {d}f_r 
\qquad f_0,f_1,\ldots,f_r\in\mnalg
\end{eqnarray}
\end{proof}


\section{An attempt at an alternative definition of $\Omega^2_\sbalg$}
\label{ch_alt_Omega2}

It seems at first that the requirement that $\sbalg$ be a generalised
algebra is unnecessarily restrictive. Especially as this is not the
case for most noncommutative versions of manifolds. One idea is to
examine the assumptions made and to see if weakening any of them would
lead to a larger choice of $\sbalg$. In this section we assume that
$\Omega^2_\sbalg$ is still bimodule of $\mnalg$, but we don't require
(\ref{2F_def_dot}). Instead we impose the condition
\begin{eqnarray}
f \theta^a \theta^b g \cdot (e_c,e_d) &=&
f g \theta^a \theta^b \cdot (e_c,e_d) =
f g Q^{ab}{}_{cd}
\qquad
\forall f,g\mnalg
\end{eqnarray}
where $Q^{ab}{}_{cd}\in\mnalg$. However, we find that if $\sbalg$ is
not a generalised algebra, then $\dim(\Omega^2_\sbalg)=0$ and we have
gained nothing. To see this we first prove.

\begin{lemma}
Given two sets of matrices $\{A^a,B^a\in M_m(\Cmpx)\}_{a=1\ldots n}$
such that
\begin{eqnarray}
\sum_{a=1}^n A^a C B^a =0 \qquad \forall C\in M_m(\Cmpx)
\label{lm_A_ind_eq}
\end{eqnarray}
Then if the set $\{A^a\}$ are independent implies all $B^a=0$.
\end{lemma}

\begin{proof}{}
Let the basis matrices be $E_{ij}$ as in lemma (\ref{1F_lm_tr})
With respect to this basis $A^a=A^a_{ij}E_{ij}$ (implicit sum on
$ij\ldots$). Putting $C=E_{ij}$ in (\ref{lm_A_ind_eq}) gives
\begin{eqnarray}
0 &=& \sum_a A^a E_{ij} B^a 
=
\sum_a A^a_{kl} E_{kl} E_{ij} E_{pq} B^a_{pq}
=
\sum_a A^a_{kl} E_{kq} \delta_{li}\delta_{jp} B^a_{pq}
\cr
&=&
\sum_a A^a_{ki} B^a_{jq} E_{kq}
\end{eqnarray}
Which since $E_{kq}$ are independent gives
\begin{eqnarray}
\sum_a A^a_{ki} B^a_{jq} = 0
\end{eqnarray}
Multiplying by $E_{ki}$ gives
\begin{eqnarray}
0 &=& 
\sum_a A^a B^a_{jq}
\end{eqnarray}
Thus implying $B^a=0$ for all $a$.
\end{proof}

\begin{proof}{of statement}
From (\ref{2F_ddf_dfdg}) then for all $f\in\mnalg$ we have
\begin{eqnarray}
0 &=& d d f = d([\lambda_a,f]\theta^a)
\cr &=&
[\lambda_b,[\lambda_a,f]]\theta^b\theta^a
- [\lambda_a,f]
\left( \theta\theta^a+
\theta^a\theta + 
\tr(\lambda^a\lambda_b\lambda_c)
\theta^b\theta^c \right)
\cr
&=&
\theta^{2} f - \theta f\theta
+ \lambda_a f \theta\theta^a
- f \lambda_a \theta\theta^a 
- f \theta^{2} + \theta f\theta
- \lambda_a f \theta\theta^a
+ f \lambda_a \theta\theta^a 
- [\lambda_a\tr(\lambda^a\lambda_b\lambda_c),f]
\theta^b\theta^c
\cr
&=&
[\theta^{2}-\eta(\lambda_b\lambda_c)\theta^b\theta^c,f]
\cr
&=&
[\eta^\perp(\lambda_a\lambda_b),f]
\theta^a\theta^b
\cr
&=&
\eta^\perp(\lambda_a\lambda_b)f
\theta^a\theta^b
\qquad\forall f\in\mnalg
\end{eqnarray}
Contracting with $(e_c,e_d)$ gives
\begin{eqnarray}
\eta^\perp(\lambda_a\lambda_b) f Q^{ab}{}_{cd} 
&=& 0
\end{eqnarray}
from the contraction of (\ref{2F_eta_theta2}) with $(e_c,e_d)$.
This is true for all $f\in\mnalg$.
Thus if $\sbalg$ is not a generalised algebra, then all
$\eta^\perp(\lambda_a\lambda_b)$ are
independent, then from the lemma above, $Q^{ab}{}_{cd}=0$ and
$\Omega^2_\mnalg$ is trivial.
\end{proof}


\section{$d$-Homomorphisms of Noncommutative algebras}
\label{ch_map}

Given any two subspaces $\sbalg,\sbalg'\subset\mnalg$ and a linear map
\begin{eqnarray}
\varphi:\sbalg\mapsto\sbalg'
\end{eqnarray}
This generates the maps
\begin{eqnarray}
\varphi_\star:\Der_\sbalg &\mapsto& \Der_{\sbalg'} \cr
\varphi_\star(\ad(h)) &=& \ad(\phi(h)) 
\qquad\forall h\in\sbalg \\
\varphi^\star:\Omega^\star_{\sbalg'} &\mapsto& \Omega^\star_\sbalg \cr
\varphi^\star(\xi')\cdot X &=& \xi\cdot(\varphi_\star X)
\qquad\forall\xi'\in\Omega^\star{\sbalg'},
X\in\Der_{\sbalg'}
\end{eqnarray}
in a similar way to that of the push forward and pull back of 
differentiable maps between manifolds.
However unlike in commutative geometry the pullback map is not in
general a $d$-homomorphisms, i.e. it does not in
general commute with the exterior derivative
$\varphi^\star d' \ne d \varphi^\star$ where
$d:\Omega^\star_{\sbalg}\mapsto\Omega^\star_{\sbalg}$ and
$d':\Omega^\star_{\sbalg'}\mapsto\Omega^\star_{\sbalg'}$.

There are some cases where they do commute.  One simple case is when
$\sbalg'$ is a subspace of $\sbalg$, if $\iota:\sbalg'\hookrightarrow
\sbalg$ then $d\iota^\star=\iota^\star d'$.  The set of relations on
products of the basis elements
$\{\lambda'_a\in\sbalg'\}_{a=1,\ldots,n'}$ making $\sbalg'$ into a
generalised algebra are, of course, a subset of the relations on
$\{\lambda_a\}$. However since in our definition of a generalised
algebra we give the possibility of ignoring some of these relationships
we cannot say that the projection matrix from (\ref{GA_P})
$P^{\prime ab}{}_{cd}$ for $\Omega^2_{\sbalg'}$ is simply the
restriction of $P^{ab}{}_{cd}$ to $\sbalg'\otimes\sbalg'$.

\subsection{Equivalent Representations}

Given $u\in GL_m(\Cmpx)$, let
\begin{eqnarray}
U:\mnalg &\mapsto& \mnalg \cr
U:\sbalg &\mapsto& \sbalg' \cr
U(h)&=& uhu^{-1}
\end{eqnarray}
This map is a bijective $d$-homomorphism, i.e. it preserves the
generalised algebraic structure.
\begin{eqnarray}
U(f+g) = U(f)+ U(g) 
\qquad\hbox{and}\qquad
U(fg) = U(f)U(g) 
\end{eqnarray}
Hence the $\alpha^{ab}_r=\alpha^{\prime ab}_r$, and $U$ may be viewed
as a map for one representation to another of the same generalised
algebra. We shall call $\sbalg$ and $\sbalg'$ {\bf equivalent}
representations. It would be nice to have some idea if given two
representations of the same generalised algebra whether they are
equivalent.  This gives rise to the following maps
\begin{eqnarray}
U_\star &:& \Der_\sbalg \mapsto \Der_{\sbalg'} \cr
U_\star(\ad h) &=& \ad(U(h)) \\
U^\star &:& \Omega^\star_{\sbalg'} \mapsto \Omega^\star_{\sbalg} \cr
(U^\star\xi')\cdot (X_1,\ldots, X_r) &=& U^{-1} (\xi\cdot (U_\star
X_1,\ldots ,U_\star X_p)) \qquad\forall \xi'\in\Omega^p_{\sbalg'},
X_1,\ldots X_p\in\Der_\sbalg
\end{eqnarray}
Note the slightly different definition of $U^\star$. 
This map has the
following properties

\begin{lemma}
If we choose the basis of $\sbalg'$ to be $\lambda'_a=U(\lambda_a)$
then $\lambda^b=u^{-1\dagger}\lambda^b u^\dagger$ and $U^\star$
preserves the co-frame $U^\star(\theta^{\prime a})=\theta^{a}$ and
$U^\star(\theta^{\prime})=\theta$. Furthermore if we choose the
$\varepsilon^{ab\ldots}$ to be equal so that
\begin{eqnarray}
\theta^{\prime a_1}\ldots\theta^{\prime a_p}
\cdot(e'_{b_1}\ldots e'_{b_p}) &=&
\theta^{a_1}\ldots\theta^{a_p}
\cdot(e_{b_1}\ldots e_{b_p}) 
\label{ERep_persv_P}
\end{eqnarray}
Then $U^\star$ preserves the exterior algebra and commutes
with the exterior derivative:
\begin{eqnarray}
U^\star(\xi'\zeta') &=& U^\star(\xi')U^\star(\zeta') 
\label{ERep_persv_prod}
\\
U^\star(d'\xi') &=& dU^\star(\xi')
\qquad\forall\xi,\zeta'\in\Omega^\star_\sbalg
\label{ERep_commut_dU}
\end{eqnarray}
\end{lemma}

\begin{proof}{}
The preservation of the co-frame is trivial. 
(\ref{ERep_persv_prod}) follows from (\ref{ERep_persv_P}).
Now
\begin{eqnarray}
U^\star [\theta',\xi'] &=&
U^\star (\theta'\xi' - (-1)^r \xi'\theta') =
U^\star(\theta')U^\star(\xi') - (-1)^r U^\star(\xi')U^\star(\theta') 
\cr
&=&
\theta U^\star(\xi') - (-1)^r U^\star(\xi')\theta =
[\theta,U^\star(\xi')]
\end{eqnarray}
and since $\langle\lambda^{\prime a},\lambda'_b\lambda'_c\rangle=
\langle\lambda^a,\lambda_b\lambda_c\rangle$ then
\begin{eqnarray}
U^\star \left(\chi(f{\theta'}^{a_1}\ldots{\theta'}^{a_p})\right)
&=&
\tfrac12 U^\star(f)\sum_{q=1}^p 
U^\star(\langle\lambda^{\prime a},\lambda'_b\lambda'_c \rangle)
U^\star({\theta'}^{a_1}\ldots{\theta'}^{a_{q-1}}
{\theta'}^b{\theta'}^c
{\theta'}^{a_{q+1}}\ldots{\theta'}^{a_p})
\cr
&=&
\tfrac12 U^\star(f)
\sum_{q=1}^p
\langle\lambda^{a},\lambda_b\lambda_c\rangle
{\theta}^{a_1}\ldots{\theta}^{a_{q-1}}
{\theta}^b{\theta}^c
{\theta}^{a_{q+1}}\ldots{\theta}^{a_p}
\cr
&=&
\chi(U^\star(f{\theta'}^{a_1}\ldots{\theta'}^{a_p}))
\end{eqnarray}
Thus from (\ref{HO_d}) we have  (\ref{ERep_commut_dU}).
\end{proof}

\subsection{The ``Lie derivative''}

As an aside we define a ``derivative''
$\Lie^\star_f:\Omega^p_\sbalg\mapsto\Omega^p_\sbalg$ for $f\in\mnalg$.
It is not obvious what r\^ole this function has. (It may be analogous
to the Lie derivative in normal commutative geometry.) All that can be
said about it is that it comes for free, that is we don't have to have
any additional structure for its definition.

Let $\sbalg'=U(\sbalg)$. As well as $U^\star$ there is another map
from $\Omega^\star_\sbalg$ to $\Omega^\star_{\sbalg'}$.  This is given
by $\phi'\phi^{-1}$ where
$\phi':\Omega^\star_u\mapsto\Omega^\star_{\sbalg'}$ is given by
(\ref{1F_comut_diag}), and $\phi^{-1}$ is a $\mnalg$-linear right
inverse of $\phi$ given by
\begin{eqnarray}
\phi^{-1}_p : \Omega^p_{\sbalg}&\mapsto&\Omega^p_u \cr
\phi^{-1}_p(f\theta^{a_1}\ldots\theta^{a_p})
&=&
P^{a_1\ldots a_p}{}_{b_1\ldots b_p}
f\theta_u^{b_1}\ldots\theta_u^{b_p}
\end{eqnarray}
Where $f\in\mnalg$, $\theta^a_u$ is given by (\ref{1F_theta_u_a})
and
\begin{eqnarray}
P^{a_1\ldots a_p}{}_{b_1\ldots b_p} &=& 
\theta^{a_1}\ldots\theta^{a_p}\cdot(e_{b_1}\ldots e_{b_p})
\end{eqnarray}
so that
\begin{eqnarray}
\phi'_p\phi^{-1}_p\theta^{a_1}\ldots\theta^{a_p}
&=&
P^{a_1\ldots a_p}{}_{b_1\ldots b_p}
\langle\lambda^{b_1},\lambda_{c_1}\rangle \ldots
\langle\lambda^{b_p},\lambda_{c_p}\rangle
\theta^{\prime c_1}\ldots\theta^{\prime c_p}
\label{LD_phiphi}
\end{eqnarray}
These two maps are not equal. Given an $f\in\mnalg$, let
$U_t(g)=e^{tf}g e^{-tf}$. We now take the derivative of the map
\begin{eqnarray}
\Omega^\star_\sbalg 
\buildrel \phi^{-1} \over \longrightarrow
\Omega^\star_u 
\buildrel \phi' \over \longrightarrow
\Omega^\star_{\sbalg'}
\buildrel U^\star_{t} \over \longrightarrow
\Omega^\star_\sbalg 
\label{LD_diag_om}
\end{eqnarray}
This is given by the ``Lie'' derivative
\begin{eqnarray}
\Lie^\star_{f} &:& \Omega^\star_\sbalg \mapsto \Omega^\star_\sbalg 
\cr
\Lie^\star_{f}(\xi) &=&
\lim_{t\to0}\tfrac1t
\left(U^\star_{t}\circ\phi'\circ\phi^{-1}(\xi)-\xi\right)
\label{LD_def_Lie}
\end{eqnarray}
which from (\ref{LD_phiphi}) is given by
\begin{eqnarray}
\Lie^\star_f(g\theta^{a_1}\ldots\theta^{a_p})
&=&
-[f,g]\theta^{a_1}\ldots\theta^{a_p}
-P^{a_1\ldots a_p}{}_{b_1\ldots b_p}
\sum_{q=1}^p \langle\lambda^{b_q},[f,\lambda_c]\rangle
\theta^{b_1}\ldots\theta^{b_{q-1}}
\theta^c\theta^{b_{q+1}}\ldots\theta^{b_q}
\end{eqnarray}
We note that
\begin{eqnarray}
U^\star(\Lie_{f}^\star \xi') &=&
\Lie_{U^\star f}^\star (U^\star\xi')
\qquad\forall f\in\mnalg,\xi'\in\Omega^\star_{\sbalg'}
\end{eqnarray}
but that in general $d\Lie^\star_f\ne\Lie^\star_f d$

\section{Examples} 
\label{ch_eg}

We give here some examples of generalised algebras which have matrix
representations and their corresponding noncommutative geometry.

\subsection{Example: The Universal Algebra}
\label{ch_eg_A0}

Let $\sbalg$ be given by $\mnalg_0\subset\mnalg$, the subspace of all
traceless matrices so that $\mnalg_0\oplus\spanone=\mnalg$ and
$n=m^2-1$.  Since $\mnalg$ is a matrix algebra then all derivations of
$\mnalg$ are in $\Der_{\mnalg_0}$. Set $P^{ab}{}_{cd}=I_{(m^2-1)^2}$
the $(m^2-1)^2\times (m^2-1)^2$ unit matrix.  In this case the map
$\phi:\Omega^\star_u\mapsto\Omega^\star_{\mnalg_0}$ given by
(\ref{1F_comut_diag}) is an isomorphism.  To see this we can choose
$\{\gamma_\mu\}_{\mu=0,\ldots,m^2-1}$ a basis for $\mnalg$ by setting
$\gamma_0=1$ and $\gamma_a=\lambda_a$ for $a=1\ldots m^2-1$ as
traceless matrices.  In this basis $\phi_1$ is given by
\begin{eqnarray}
\phi_1 \left(
\sum_{\mu,\nu=0}^{m^2-1} \xi_{\mu\nu}
\gamma_\mu\otimes\gamma_\nu
\right)=
\sum_{a,b=1}^{m^2-1} \xi_{ab}
\lambda_ad\lambda_b
+
\sum_{b=1}^{m^2-1} \xi_{0b}
d\lambda_b
\end{eqnarray}
The inverse of this map can be calculated since
\begin{eqnarray}
-\sum_{a,b=1}^{m^2-1} \xi_{ab}\lambda_a\lambda_b
-\sum_{b=1}^{m^2-1} \xi_{0b}\lambda_b
=
\sum_{a=1}^{m^2-1} \xi_{a0}\lambda_a
+ \xi_{00}
\end{eqnarray}
This is extended for all $\phi_p$.  The space of $p$-forms is now a
free bimodule over $\mnalg$ of rank $(m^2-1)^p$, and all the co-frame
basis elements $\theta^{a_1}\ldots\theta^{a_p}$ are independent. The
1-form $\theta$ is given by $\theta_u$ in (\ref{1F_theta_u}).

We can now view any other noncommutative geometry given by the
subspace $\sbalg\subset\mnalg_0$ as being a sub noncommutative
geometry as stated in section \ref{ch_map}. The maps
$\iota:\sbalg\mapsto\mnalg_0$ and $\eta:\mnalg_0\mapsto\sbalg$ induce
the pullbacks $\iota^\star:\Omega^\star_{\mnalg_0}\mapsto
\Omega^\star_\sbalg$ and $\eta^\star:\Omega^\star_\sbalg\mapsto
\Omega^\star_{\mnalg_0}$. By identifying $\Omega^\star_u$ and
$\Omega^\star_{\mnalg_0}$ then
$\iota^\star=\phi_\sbalg:\Omega^\star_u\mapsto\Omega^\star_\sbalg$,
and so, of course, commutes with $d$. It is easy to show that
$\eta^\star(\theta^a)=\theta^a_u$ given by (\ref{1F_theta_u_a}).

We can also view \cite[chapter 3]{Madore_book} $\mnalg_0$ as the
fundamental representation of the Lie algebra $sl(m)$. For this we
must choose the elements $\beta^r_{ab}$ so that
(\ref{eglie_theta_anticom}) below holds.

\subsection{Example: The Lie Algebra}
\label{ch_lie}

A standard example of a generalised algebra is the case of a Lie
algebra. This case has been studied in detail \cite{Madore_book}, 
especially when $\sbalg$ is a representation of
$su(2)$, which has been shown to be a non-commutative approximation to
the sphere and $su(4)$ which is an analogue of the Euclidianised
compactified Minkowski space.

If $\sbalg$ is a representation of a Lie group of dimension $n$ 
then 
\begin{eqnarray}
[\lambda_a,\lambda_b] = C^c{}_{ab}\lambda_c \in\sbalg
\label{Lie_fund_eqn}
\end{eqnarray}
for $a,b=1\ldots n$, where $C^c{}_{cd}$ are the structure constants.
This make a total of $\tfrac12n(n-1)$ independent
equations. There are also the Casimir operators 
\begin{eqnarray}
\sum_{a=1}^{n'} \lambda'_i \lambda'_i \in\spanone
\end{eqnarray}
for any orthogonal basis $\{\lambda'_i\}_{i=1\ldots n'}$ of either
$\sbalg$ or any sub Lie algebra $\sbalg'\subset\sbalg$.  It is usual
to ignore all these equations, and take only those given by
(\ref{Lie_fund_eqn}).  Thus the rank of the generalised algebra
$\rank=\tfrac12n(n-1)$.  Hence it is easier to replace $r$ by the pair
$(c,d)$ with $c<d$. Thus
\begin{eqnarray}
\alpha^{ab}_{cd} &=& \delta^a_c\delta^b_d - \delta^a_d\delta^b_c
\cr 
\beta^{cd}_{ef} &=& \tfrac12
(\delta^c_e\delta^d_f - \delta^c_f\delta^d_e )
\qquad\hbox{for }c<d
\end{eqnarray}
so
\begin{eqnarray}
\theta^a\theta^b+\theta^b\theta^a &=& 0
\label{eglie_theta_anticom}
\end{eqnarray}
We also note that ${h^\dagger}\in\sbalg$
for all $h\in\sbalg$. Thus we can choose
$\lambda_a$ to be Hermitian, or anti-Hermitian. We get the same
results if the $\lambda_a$'s mutually commuted. We would then set
$C^a{}_{bc}=0$ in (\ref{Lie_fund_eqn}).


\subsection{Example: $q$-Deformed Algebra }
\label{ch_eg_qyx}

A $q$-deformed algebra $\mnalg$ is generated by the elements
$x,y\in\mnalg$ where $xy=qyx$.  We can find a $M_m(\Cmpx)$
representation for a $q$-deformed algebra if and only if $q^m=1$. In
order that $q\to0$ as $m\to\infty$ let $q=e^{2\pi i/m}$.  A
representation is then given by
\begin{eqnarray}
x=
\left(
\begin{array}{c|c}
\strut
\;0\; & \;I_{m-1} \;\\
\hline
\strut
1 & 0
\end{array}
\right)\ ,\qquad
y =
\hbox{diag}(1,q,q^2,\ldots,q^{m-1})
\end{eqnarray}
where $I_{m-1}\in M_{m-1}(\Cmpx)$ is the identity matrix.
We see that $x$ and $y$ are non degenerate, traceless matrices.
Since $x^\dagger x=y^\dagger y=1$ let
\begin{eqnarray}
\lambda_1 = x,\quad
\lambda_2 = y,\quad
\lambda^1 = \tfrac1m x,\quad
\lambda^2 = \tfrac1m y
\end{eqnarray}
Explicit forms of $\theta^1,\theta^2$ are given in
\cite{Madore_Dim1}.
We note that $\tr(x^ay^b)=0$ if either $a$ or $b$ is not a multiple of
$m$. 
Thus we have:
\begin{eqnarray}
\langle \lambda^a, \lambda_b\lambda_c \rangle = 0,
\hbox{\qquad and \qquad}
\tr(\lambda_a\lambda_b)=0
\qquad\forall a,b,c=1,2
\end{eqnarray}
Equations (\ref{2F_dtheta_a}) and (\ref{2F_dtheta_theta_sqr}) become
\begin{eqnarray}
&&d\theta^a=-[\theta,\theta^a] \\
&&d\theta=-\theta^2
\end{eqnarray}
The space of 2-forms $\Omega^2_\sbalg$ is given by
(\ref{2F_eta_theta2})
\begin{eqnarray}
\lambda_1\lambda_1\theta^{1}\theta^{1} +
\lambda_2\lambda_1
(q\theta^{1}\theta^{2} + \theta^{2}\theta^{1}) +
\lambda_2\lambda_2\theta^{2}\theta^{2} &=& 0
\end{eqnarray}
Which implies that the rank of $\Omega^1_\sbalg=2$ and
the rank of $\Omega^1_\sbalg=1$ with
\begin{eqnarray}
\theta^{1}\theta^{1}=\theta^{2}\theta^{2}=0,\qquad
q\theta^{1}\theta^{2} + \theta^{2}\theta^{1} =0
\end{eqnarray}


\subsection{Example: The ``Fuzzy Ellipsoid'' }

In the previous three examples the generalised algebra and associated
non commutative geometry were already well established.  Here we give
a simple generalised algebra of rank 1 which has not been studied
before. We have called it the ``Fuzzy Ellipsoid'' since it is based on
the fuzzy sphere with two of the three elements of $\sbalg$ unchanged.

Let $\{J_1,J_2,J_3\}$ be an $M_m(\Cmpx)$ Hermitian representation of
$su(2)$ such 
that $[J_i,J_j]=i\varepsilon_{ijk}J_k$. 
Let 
\begin{eqnarray}
\sbalg &=& \span\{\lambda_1,\lambda_2,\lambda_3\}
\end{eqnarray}
where $\lambda_1 = -i\kappa J_1,\ \lambda_2 = -i\kappa J_2,$
and 
\begin{eqnarray}
i\lambda_3 &=& 
\alpha^{11} \lambda_1\lambda_1 +
\alpha^{12} \lambda_1\lambda_2 +
\alpha^{21} \lambda_2\lambda_1 +
\alpha^{22} \lambda_2\lambda_2
-
\tfrac1{12} \kappa^2 m(m^2-1) (\alpha^{11} - \alpha^{22})
\end{eqnarray}

In this space we have $\dim_\mnalg(\Omega^1_\sbalg)=3$ and
 $\dim_\mnalg(\Omega^2_\sbalg)=1$ consisting of the span of the
element
\begin{eqnarray}
\frac{\theta^{1}\theta^{1}}{\alpha^{11}} =
\frac{\theta^{1}\theta^{2}}{\alpha^{12}} =
\frac{\theta^{2}\theta^{1}}{\alpha^{21}} =
\frac{\theta^{2}\theta^{2}}{\alpha^{22}} 
\end{eqnarray}
with $\theta^{a}\theta^{b}=0$ otherwise.

The elements $\tr(\lambda_a\lambda_b)$ and
$\langle\lambda^a,\lambda_b\lambda_c\rangle$ can be calculated.
However these are simply long 6th order multipolynomials in $m$ and
$\alpha^{ab}$ which don't give any information.

\section{Discussion}

It would be nice to know which of the results discussed in this paper
can be generalised to infinite dimensional algebras, or alternatively
to infinite dimensional representation of generalised algebras. This
is necessary for example in the $q$-deformed algebras when there
doesn't exist an $m\in\Intg$ such that $q^m=1$, and for any
representation of the Heisenberg algebra $[p,x]=i$.  From the counter
examples on page \pageref{1F_counter_eg} one cannot assume that
$\Omega^1_\sbalg$ has a co-frame.  Also since the trace of an operator
is not in general defined this will cause further problems as
lemma \ref{1F_lm_tr} cannot be applied.  We know from
\cite{Madore_Dim1,Madore_Dim2} that there will still exist
restrictions on $\sbalg$ equivalent to demanding that it is a
generalised algebra.

There has been much work recently on the concept of linear connections
and curvature in noncommutative geometry
\cite{Madore_Dim1,Madore_Dim2,DubVio1,DubVio2,DubVio3,Mourad} and
\cite[chapter 3.5]{Madore_book} 
This work has been limited in the main to established algebras, such
as Lie and $q$-deformed and quantum algebra. It
would be nice to extend this work for Generalised Algebras.

As stated in section \ref{ch_map} it would be nice to have some
theorems (in line with those for Lie algebras) that ascertain when a
generalised algebra has a matrix representation, and if given two such
representation when they are equivalent.

\vskip 2em

{\bf\LARGE Acknowledgements}

The author would like to thank John Madore and Jihad Mourad for useful
discussions which motivated this work.  The author would also like to
thank the Royal Society of London for a European Junior Fellowship,
and Richard Kerner and the Laboratoire de Gravitation et Cosmologie
Relativistes, Paris~VI for their hospitality.

\end{document}